# Optimization of drug controlled release from multi-laminated devices based on the modified Tikhonov regularization method

Xinming Zhang[*], Ling Ma

(School of Science, Harbin Institute of Technology (Shenzhen), Shenzhen 518055, China)

**Abstract:** From the viewpoint of inverse problem, the optimization of drug release based on the multi-laminated drug controlled release devices has been regarded as the solution problem of the diffusion equation initial value inverse problem. In view of the ill-posedness of the corresponding inverse problem, a modified Tikhonov regularization method is proposed by constructing a new regularizing filter function based on the singular value theory of compact operator. The convergence and the optimal asymptotic order of the regularized solution are obtained. Then the classical Tikhonov regularization method and the modified Tikhonov regularization method are applied to the optimization problem of the initial drug concentration distribution. For three various desired release profiles (constant release, linear decrease release and linear increase followed by a constant release profiles), better results can be obtained by using the modified Tikhonov regularization method. The numerical results demonstrate that the modified Tikhonov regularization method not only has the optimal asymptotic order, but also is suitable for the optimization and design of multi-laminated drug controlled release devices.

## 1. Introduction

It is well known that the controlled release device is usually used to regulate the release of active material, in order to maintain a predetermined concentration of the active material for a specified period of time. In recent years, the technique has been widely used in many fields, including drug, food and cosmetics et cl [1,2] due to its safety and effectiveness. Especially, in drug research field, the controlled release device has been of interest to many scholars. Generally, for drug release behavior, the burst effect is not an desired delivery profile which will cause overdose of drug and can result in systemic toxicity [3]. To overcome this undesired behavior, various matrix geometries have been employed over the last three decades. Among them, an effective and simple device geometry is the multi-layer matrix devices [4-6] where the matrix core, containing the drug, is covered by one or more modulating layers that act as rate-controlling barriers. Researches have revealed that for the multi-layer matrix devices, the reasonable initial drug loading can efficiently control the burst effect [7].

For understanding and simulating the real release mechanisms from these complicated multi-laminated matrix devices better, application of mathematical models is a good choice [8-10]. In 1961, Takeru Higuchi[11], who is regarded as the "father of mathematical modeling of drug release", proposed a simple but surprising equation to quantify the drug release from monolithic dispersions with slab geometry: the famous Higuchi equation. After that, more and more scientists paid their attentions to the mathematical modeling of drug release, which can significantly contribute to the product development and mechanism understanding. Especially, over the last thirty years, there are many new progresses for different controlled release devices. Helbling et al[12] proposed an analytical solutions for the case of controlled dispersed-drug release from planar non-erodible polymeric matrices based on refined integral method in 2010, and then presented a mathematical modeling of controlled release of drug from torus-shaped single-layer devices and derived analytical solutions based on the pseudo-steady state approximation in 2011. Streubel[13] developed a new multi-layer matrix tablets to achieve bimodal drug release profiles and investigated the process and formulation parameters affecting the resulting release rates with theophylline and acetaminophen as model drugs. Yin and Li[14] introduced the fractional calculus to model the anomalous diffusions and presented some new mathematical models for the release of drug from both non-degradable and degradable slab matrices. Wu and Zhou[15] used the finite element method to study the effect of several factors on the kinetics of diffusional drug release from complex geometries. Nevertheless, most of the previous methods focus on how to predict release profiles based on the known initial parameters. Rather less effort has been spent on identifying the optimal values of the available control parameters, in multi-layer devices, to make the drug release profile as close to the desired release profile as possible. As far as we known, in 1998 and 1999, articles [16, 17] made relatively big contribution to this area firstly. In these two manuscripts, Lu and his co-workers employed a formal optimization approach to correctly determine the initial drug concentration in the layer so as to coincide the required release profile as much as possible. Georgiadis and Kostoglou[18] presented a systematic optimization framework to achieve desired release rates in drug delivery devices using multi-laminated layers based on a simple mathematical model. Nauman et al[19] proposed a readily manufacturable design that consists of two or three layers having different initial loadings of the chemical. By adjusting the parameters in this device, a variety of release profiles, i.e. flux versus time, can be achieved. In these previous articles, the optimizations of the available control parameters have been performed in the frame of optimization. However, we can also treat this problem from another viewpoint, that is, the viewpoint of inverse problem. In fact, many practical problems in different fields can be reduced into the solution of inversion problems, such as geophysical exploration, medical imaging and so on. In this article, from the viewpoint of the inverse problem, the optimization problem of drug release based on the multi-laminated drug controlled release devices has been transformed into the solution problem of diffusion equation initial value inverse problem. And a classical regularization method, that is Tikhonov regularization method[20], and its variant have been attempted to solve it. In 1963, Tikhonov proposed the Tikhonov regularization method firstly. After that, many scholars paid their attentions on this rather effective method and developed many different variants. A. Neubauer [21]studied Tikhonov regularization as a stable method for approximating the solutions of non-linear ill-posed problems in detail; Fu et.al[22] applied the Fourier method to solve some ill-posed problems and systematically considered a posteriori choice of the regularization parameter; Zheng and Wei[23] used a spectral regularization method to solve the Cauchy problem of TFADE based on the solution given by the Fourier method; Zhao[24] introduced a mollification method based on expanded Hermit functions to solve the ill-posedness of the problem; Chen, Fu and Feng[25] presented an optimal filtering method to approximate a Cauchy problem for the Helmholtz equation in a rectangle and showed the Hölder type error estimate; Bonesky et.al[26] applied an adaptive wavelet algorithms to solve an inverse parabolic problem that stems from the industrial process of melting iron ore in a steel furnace; Jin[27] formulated an inexact Newton–Landweber iteration method to solve nonlinear inverse problems in Banach spaces by making use of duality mappings; Bauer et al[28] studied the convergence of regularized Newton methods applied to nonlinear operator equations in Hilbert spaces if the data are perturbed by random noise.

In this article, a modified Tikhonov regularization method is proposed by constructing a new regularizing filter function based on the singular value theory of compact operator. The convergence and the optimal asymptotic order of the regularized solution are obtained. Then the Tikhonov regularization method and the modified Tikhonov regularization method are applied to the optimization of the initial drug concentration distribution of the multi-laminated drug controlled release devices. In the reminder of the paper, it is organized as follows. In Section 2, we describe the mathematics model for the multi-laminated drug controlled release devices and the corresponding inverse problem frame for the optimization of drug controlled release. In Section 3, a modified Tikhonov regularization method is described in detail. This is followed by numerical simulation in Section 4 and Section 5, lastly the conclusion is indicated in Section 6.

## 2. Mathematics Model and inverse problem frame

### 2.1 Mathematics Modeling

Fig.1 shows a multi-laminated drug controlled release device, which has N layers. The device has a thickness L and initial drug concentration $V(X)$. It makes contact with the environment through the rightmost layer and is sealed at the other side by an impermeable layer. It is assumed that the device is not significantly swelling and eroding during drug release. Here, $C_i \ (i=1,2,\cdots,N)$ are the drug concentrations of each layer, respectively. The present analysis focuses only on the case of low drug concentration.

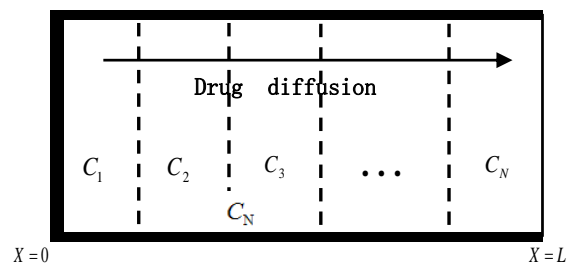

Fig.1 Drug release from a multi-laminated drug controlled release device

Mathematically, this problem can be modeled as one dimensional partial differential equation according to Fick's second law of diffusion:

$$\frac{\partial C}{\partial \tau} = \frac{\partial}{\partial X}\left(D\frac{\partial C}{\partial X}\right) \qquad (1)$$

where C is the drug concentration, $\tau$ is the release time and X is the position for one-dimensional diffusional processes, D is the diffusivity.

We assume zero flux at the interface with impermeable layer and zero concentration at the interface with the environment. The boundary conditions are therefore

$$\left.\frac{\partial C}{\partial X}\right|_{X=0} = 0, \quad \tau > 0, \qquad (2)$$

$$C(\tau, L) = 0, \quad \tau > 0, \qquad (3)$$

The initial condition is imposed:

$$C(0, X) = V(X), \quad \tau = 0, \ 0 < X < L, \qquad (4)$$

The flux of drug is also defined as:

$$J(\tau, L) = -D_{(X=L)}\left.\frac{\partial C}{\partial X}\right|_{X=L}, \quad \tau > 0, \qquad (5)$$

**2.2 The Inverse Problem Frame**

For the mathematical model (Eqs.(1)-(5)), if the initial conditions and the boundary conditions are known, the processing to compute the concentration distribute functions $C(\tau, X)$ is a forward problem, which is well-posed problem. Conversely, how to identify the initial conditions (4) is a classical inverse problem, if we know the boundary conditions (2)-(3) and additional condition (5). Generally speaking, the inverse problem is ill-posed and always be solved with some special algorithms, such as the Tikhonov regularization method.

Assume that the diffusivity is constant, we use the following dimensionalization:

$c = C/C_0$, $v = V(X)/C_0$, $x = X/L$, $t = D_0\tau/L^2$, $j = JL/D_0 C_0$, $d = D/D_0 = 1$.

where $C_0$ is a reference concentration and $D_0$ is a reference diffusivity.

The mathematical model can be transformed into the following:

$$\frac{\partial c}{\partial t} = \frac{\partial^2 c}{\partial x^2} \qquad (6)$$

Boundary conditions:

$$\left.\frac{\partial c}{\partial x}\right|_{x=0} = 0, \ t > 0, \qquad (7)$$

$$c(t,1) = 0, \ t > 0, \qquad (8)$$

Initial conditions: $\quad c(0, x) = v(x), \ t = 0, \ 0 < x < 1, \qquad (9)$

Additional conditions:

$$j(t,1) = -\left.\frac{\partial c}{\partial x}\right|_{x=1} = j(t), \quad t > 0, \tag{10}$$

In fact, the problem to determine the initial condition $v(t)$ based on the previous mathematical modeling is a diffusion equation inverse problem, and can come down to a solution of Fredholm integral equation of first kind.

The first step is to solve the equations (6)-(9), which is diffusion equation initial boundary value problem. The method of separating variables leads to the analytical solution:

$$c(t,x) = \sum_{k=0}^{\infty} 2e^{-\left[\left(k+\frac{1}{2}\right)\pi\right]^2 t} \cos\left(k+\frac{1}{2}\right)\pi x \int_0^1 v(x)\cos\left(k+\frac{1}{2}\right)\pi x\, dx \tag{11}$$

The flux $j(t)$ in additional conditions (10) is determined by differentiation of Equaion (11) with respect to $x$

$$\begin{aligned} j(t) = -\left.\frac{\partial c}{\partial x}\right|_{x=1} &= \sum_{k=0}^{\infty} 2e^{-\left[\left(k+\frac{1}{2}\right)\pi\right]^2 t}\left(k+\frac{1}{2}\right)\pi \sin\left(k+\frac{1}{2}\right)\pi x \int_0^1 v(x)\cos\left(k+\frac{1}{2}\right)\pi x\, dx \\ &= \sum_{k=0}^{\infty} 2(-1)^k \left(k+\frac{1}{2}\right)\pi e^{-\left[\left(k+\frac{1}{2}\right)\pi\right]^2 t} \int_0^1 v(x)\cos\left(k+\frac{1}{2}\right)\pi x\, dx \\ &= 2\int_0^1 \left[\sum_{k=0}^{\infty}(-1)^k\left(k+\frac{1}{2}\right)\pi e^{-\left[\left(k+\frac{1}{2}\right)\pi\right]^2 t}\cos\left(k+\frac{1}{2}\right)\pi x\right] v(x)\, dx \end{aligned} \tag{12}$$

Using the equation (12), the previous inverse problem for determining initial condition can be transformed into the following Fredholm integral equation of first kind:

$$2\int_0^1 \left[\sum_{k=0}^{\infty}(-1)^k\left(k+\frac{1}{2}\right)\pi e^{-\left[\left(k+\frac{1}{2}\right)\pi\right]^2 t}\cos\left(k+\frac{1}{2}\right)\pi x\right] v(x)\,dx = j(t), \tag{13}$$

where $v(x)$ is unknown function and the kernel function is

$$k(x,t) = \sum_{k=0}^{\infty}(-1)^k\left(k+\frac{1}{2}\right)\pi e^{-\left[\left(k+\frac{1}{2}\right)\pi\right]^2 t}\cos\left(k+\frac{1}{2}\right)\pi x \tag{14}$$

### 3. The modified Tikhonov regularization method

Based on the concept of the regularizing filter function proposed by Kirsch [29], the filter function corresponding to the Tikhonov regularization method is: $q(\alpha,\mu) = \dfrac{\mu^2}{\alpha + \mu^2}$. In 2006, Li[30] proposed an improved Tikhonov regularization method based on a modified regularizing filter function: $q(\alpha,\mu) = \dfrac{\mu^\sigma}{(\alpha + \mu^{\sigma r})^{\frac{1}{r}}}, r > 0, \sigma \geq 1,$ and proved that the regularization solution can achieve the optimal asymptotic convergence rate by selecting reasonable regularization parameters. In this paper, a new regularization Tikhonov method is presented by combining the previous modified regularization method with the truncated singular value decomposition (TSVD) regularization method. The convergence and the optimal asymptotic order of the new regularized solution are also obtained.

As we all know, we have the following lemma for the regularizing filter function:

**Lemma 1.**[29] Let X, Y are both Hilbert Spaces, $K: X \to Y$ be compact with singular system $(\mu_i, x_i, y_i)$ and $q:(0,+\infty) \times (0, \|K\|] \to R$ be a function with the following properties:

（1）$|q(\alpha,\mu)| \leq 1, \forall \alpha \in (0,+\infty), \forall \mu \in (0, \|K\|]$;

（2）For $\forall \alpha \in (0,+\infty)$, there exists $c(\alpha) > 0,$ such that $|q(\alpha,\mu)| \leq c(\alpha)\mu, \forall \mu \in (0, \|K\|]$;

（3）$\lim\limits_{\alpha \to 0} q(\alpha,\mu) = 1$, $\forall \mu \in (0, \|K\|]$.

Then the operator $R_\alpha: Y \to X$, defined by

$$R_\alpha y = \sum_{i=1}^{\infty} \frac{q(\alpha, \mu_i)}{\mu_i}(y, y_i)x_i, \tag{15}$$

is a regularization strategy with $\|R_\alpha\| \leq c(\alpha)$. The function $q(\alpha,\mu)$ with the above properties is called regularizing filer function for K.

For the famous Tikhonov regularization method, the Tikhonov regularized solution is the minimum value of the Tikhonov functional and can be expressed as the following formula by applying the singular system of the compact operator:

$$x_\alpha^\delta = R_\alpha y_\delta = \sum_{i=1}^{\infty} \frac{q(\alpha, \mu_i)}{\mu_i}(y_\delta, y_i)x_i, \tag{16}$$

where $q(\alpha,\mu) = \dfrac{\mu^2}{\alpha + \mu^2}, \alpha > 0, 0 < \mu \leq \|K\|$ is the Tikhonov regularization filter function.

Similarly, we can construct a novel regularization method by giving a new regularizing filter function.

**3.1 The construction of a new regularizing filter function**

The filter function corresponding to the modified regularization method proposed by literature [30] is:

$$q(\alpha,\mu) = \frac{\mu^\sigma}{(\alpha+\mu^{\sigma r})^{\frac{1}{r}}}, r>0, \sigma \geq 1, \tag{17}$$

The filter function of the truncated singular value decomposition regularization method is:

$$q(\alpha,\mu) = \begin{cases} 1 & \mu \geq \alpha \\ 0 & \mu < \alpha \end{cases}. \tag{18}$$

A new modified regularizing filter function can be constructed as follows by combining the previous two approach:

$$q(\alpha,\mu) = \begin{cases} 1 & \mu^{\sigma r} \geq \alpha \\ \dfrac{\mu^\sigma}{(\alpha+\mu^{\sigma r})^{\frac{1}{r}}} & \mu^{\sigma r} < \alpha \end{cases} \tag{19}$$

where $\alpha > 0, 0 < \mu \leq \|K\|, r > 0, \sigma \geq 1$.

Based on formula (19), the large singular value cannot be modified, meanwhile the small singular can also be filtered. Furthermore, the accuracy of the solution will not be affected. It can be proved that:

**Theorem 1. The function $q(\alpha,\mu)$ defined by (19) is a regularizing filter function that lead to optimal regularization strategy.**

**Proof.**

1) For $\mu^{\sigma r} \geq \alpha$ one has $q(\alpha,\mu)=1 \leq 1$; for $\mu^{\sigma r} < \alpha$, since $\mu^\sigma = (\mu^{\sigma r})^{1/r} < (\alpha+\mu^{\sigma r})^{1/r}$, so $q(\alpha,\mu) = \dfrac{\mu^\sigma}{(\alpha+\mu^{\sigma r})^{1/r}} < 1$.

2) For $\mu^{\sigma r} \geq \alpha$, one has $\dfrac{\mu^{\sigma r}}{\alpha} \geq 1$, that is $\dfrac{\mu}{\alpha^{1/\sigma r}} \geq 1$, so $|q(\alpha,\mu)| = 1 \leq \dfrac{1}{(\alpha)^{1/\sigma r}} \mu$

For $\mu^{\sigma r} < \alpha$, $q(\alpha,\mu) = \dfrac{\mu^\sigma}{(\alpha+\mu^{\sigma r})^{1/r}}$. The choice $\sigma = 1$ leads to

$$q(\alpha,\mu) = \frac{\mu}{(\alpha+\mu^r)^{1/r}} \leq \frac{\mu}{\alpha^{1/r}} \tag{20}$$

For $\sigma > 1$, let $p = \sigma, q = \dfrac{\sigma}{\sigma-1}$. Thanks to the Yong inequality $\alpha + \mu^{\sigma r} \geq \alpha^{1/p} \cdot \mu^{\sigma r/q}$, we have $\alpha + \mu^{\sigma r} \geq \alpha^{1/\sigma} \cdot \mu^{r(\sigma-1)}$, which yields

$$\frac{\mu^\sigma}{(\alpha+\mu^{\sigma r})^{1/r}} \leq \frac{\mu^\sigma}{\alpha^{1/\sigma r} \cdot \mu^{\sigma-1}} = \frac{\mu}{\alpha^{1/\sigma r}} \tag{21}$$

Thus the inequality $q(\alpha, \mu) \leq \dfrac{\mu}{\alpha^{1/\sigma r}}$ holds for $\forall \sigma \geq 1$. That is the property (2) in Lemma holds with $c(\alpha) = \dfrac{1}{\alpha^{1/\sigma r}}$, for $\forall \alpha > 0$

3) It is obvious that $q(\alpha, \mu) = 1$, as $\alpha \to 0$.

From **Lemma 1**, we know that the function $q(\alpha, \mu) = \begin{cases} 1 & \mu^{\sigma r} \geq \alpha \\ \dfrac{\mu^{\sigma}}{(\alpha + \mu^{\sigma r})^{\frac{1}{r}}} & \mu^{\sigma r} < \alpha \end{cases}$ is a

regularizing filter function, and the regularization operator is $R_\alpha : Y \to X$

$$R_\alpha y = \sum_{i=1}^{\infty} \dfrac{q(\alpha, \mu_i)}{\mu_i}(y, y_i) x_i \tag{22}$$

**3.2 Error analysis of the regularized solutions**

**Theorem 2.** Let $x^+ = (K^*K)^v z \in R(K^*K)^v$, $z \in X$, with $\|z\| \leq E$, for the regularization operators $R_\alpha : Y \to X$ from (22), we choose the regularized parameters $\alpha(\delta) = c(\dfrac{\delta}{E})^{\frac{\sigma r}{2v+1}}$ for some $c > 0$, then the following estimate holds:

$$\|x^\delta_{\alpha(\delta)} - x^+\| = O(\delta^{\frac{2v}{2v+1}}). \tag{23}$$

where $x^+$ is the solution of $Kx = y$, $x^\delta_{\alpha(\delta)} := R_\alpha y^\delta$ is the approximation of $x$, $K^*$ denotes the adjoint of $K$.

**Proof:** The error between $x$ and $x^\delta_{\alpha(\delta)}$ is

$$\|x^\delta_\alpha - x^+\| \leq \|R_\alpha\| \cdot \delta + \|R_\alpha y - x^+\|, \tag{24}$$

**Theorem 1** yields $\|R_\alpha\| \leq c(\alpha) = \dfrac{1}{\alpha^{1/\sigma r}}$. From $R_\alpha Kx = \sum_{i=1}^{\infty} \dfrac{q(\alpha, \mu_i)}{\mu_i}(Kx, y_i)x_i$, $x = \sum_{i=1}^{\infty}(x, x_i)x_i$

and $(Kx, y_i) = (x, K^* y_i) = \mu_i(x, x_i)$, we conclude that

$$\begin{aligned}
\|R_\alpha y - x^+\|^2 &= \sum_{i=1}^{\infty} |q(\alpha,\mu_i)-1|^2 \cdot |(x^+, x_i)|^2 \\
&= \sum_{i=1}^{\infty} |q(\alpha,\mu_i)-1|^2 \cdot |((K^*K)^v z, x_i)|^2 \\
&= \sum_{i=1}^{\infty} |q(\alpha,\mu_i)-1|^2 \cdot \left|\left(\sum_{j=1}^{\infty} \mu_j^{2v}(z, x_j) x_j, x_i\right)\right|^2 \\
&= \sum_{i=1}^{\infty} |q(\alpha,\mu_i)-1|^2 \cdot |\mu_i^{2v}(z, x_i)|^2 \\
&= \sum_{i=1}^{\infty} |q(\alpha,\mu_i)-1|^2 \cdot \mu_i^{4v} \cdot |(z, x_i)|^2
\end{aligned} \tag{25}$$

For $\mu^{\sigma r} < \alpha$, $|q(\alpha,\mu_i)-1| < 1$ because $0 < q(\alpha,\mu) < 1$, thus

$$|q(\alpha,\mu_i)-1| \cdot \mu_i^{2v} < \mu_i^{2v} = (\mu_i^{\sigma r})^{2v/\sigma r} < \alpha^{2v/\sigma r}, \tag{26}$$

For $\mu^{\sigma r} \geq \alpha$, $q(\alpha,\mu) = 1$, so

$$|q(\alpha,\mu_i)-1| \cdot \mu_i^{2v} = 0 \cdot \mu_i^{2v} = 0 < \alpha^{2v/\sigma r}, \tag{27}$$

Then $|q(\alpha,\mu_i)-1| \cdot \mu_i^{2v} < \alpha^{2v/\sigma r}$ holds for all two cases.

Therefore

$$\|R_\alpha y - x^+\|^2 < \alpha^{4v/\sigma r} \sum_{i=1}^{\infty} |(z, x_i)|^2 = \alpha^{4v/\sigma r} \cdot \|z\|^2 \leq \alpha^{4v/\sigma r} \cdot E^2, \tag{28}$$

that is $\|R_\alpha y - x^+\| < \alpha^{2v/\sigma r} \cdot E$. Thus we have shown that

$$\|x_\alpha^\delta - x^+\| \leq \frac{1}{\alpha^{1/\sigma r}} \cdot \delta + \alpha^{\frac{2v}{\sigma r}} \cdot E, \tag{29}$$

The choice of $\alpha(\delta) = c(\frac{\delta}{E})^{\frac{\sigma r}{2v+1}}$ leads to the corresponding estimate

$$\begin{aligned}
\|x_{\alpha(\delta)}^\delta - x^+\| &\leq \delta \cdot [c(\frac{\delta}{E})^{\frac{\sigma r}{2v+1}}]^{-\frac{1}{\sigma r}} + [c(\frac{\delta}{E})^{\frac{\sigma r}{2v+1}}]^{\frac{2v}{\sigma r}} \cdot E \\
&= \delta \cdot [c^{-\frac{1}{\sigma r}} \cdot (\frac{\delta}{E})^{-\frac{1}{2v+1}}] + c^{\frac{2v}{\sigma r}} \cdot (\frac{\delta}{E})^{\frac{2v}{2v+1}} \cdot E \\
&= c^{-\frac{1}{\sigma r}} \cdot E^{\frac{1}{2v+1}} \cdot \delta^{1-\frac{1}{2v+1}} + c^{\frac{2v}{\sigma r}} \cdot \delta^{\frac{2v}{2v+1}} \cdot E^{1-\frac{2v}{2v+1}} \\
&= (c^{-\frac{1}{\sigma r}} + c^{\frac{2v}{\sigma r}}) \cdot E^{\frac{1}{2v+1}} \cdot \delta^{\frac{2v}{2v+1}}
\end{aligned} \tag{30}$$

This yields

$$\|x_{\alpha(\delta)}^\delta - x^+\| = O(\delta^{\frac{2v}{2v+1}}), \tag{31}$$

## 4. Solution of the integral equations of first kind

In this section, we will give two examples to show the effectiveness of our method.

**Example 4.1** The integral equation of first kind

$$\int_0^1 (1+ts)e^{ts}x(s)\mathrm{d}s = y(t), \quad 0 \leq t \leq 1 \tag{32}$$

with the right-hand side and the solution given by $y(t) = e^t$ and $x(t) = 1$.

We discretize the integral equation by the compound trapezoidal formula and obtain linear system

$$\mathbf{AX} = \mathbf{b} \tag{33}$$

where $\mathbf{A} \in R^{100 \times 100}$ and the error-free right hand side $\mathbf{b} = \left[ y(t_0), \cdots, y(t_m) \right]^T \in R^{100}$. The associated contaminated vector $\tilde{\mathbf{b}}$ is given by the following formula

$$\tilde{\mathbf{b}} = \left[ y(t_0), \cdots, y(t_m) \right]^T + \delta \times rand(m+1, 1) \tag{34}$$

where $\delta = 0.001$.

We solve this equation with the Tikhonov regularization method and the modified Tikhonov regularization method, respectively. The regularization parameter is chosen by using the L-curve method, and the corresponding curve is shown in Fig.2 and Fig.4. Fig. 3 depicts the comparison between the results obtained by the classical Tikhonov regularization method and exact solution, meanwhile Fig.5 gives the comparison between the results obtained by the modified Tikhonov regularization method and exact solution. As we can see that the modified Tikhonov regularization method are more effective than the classical Tikhonov regularization method.

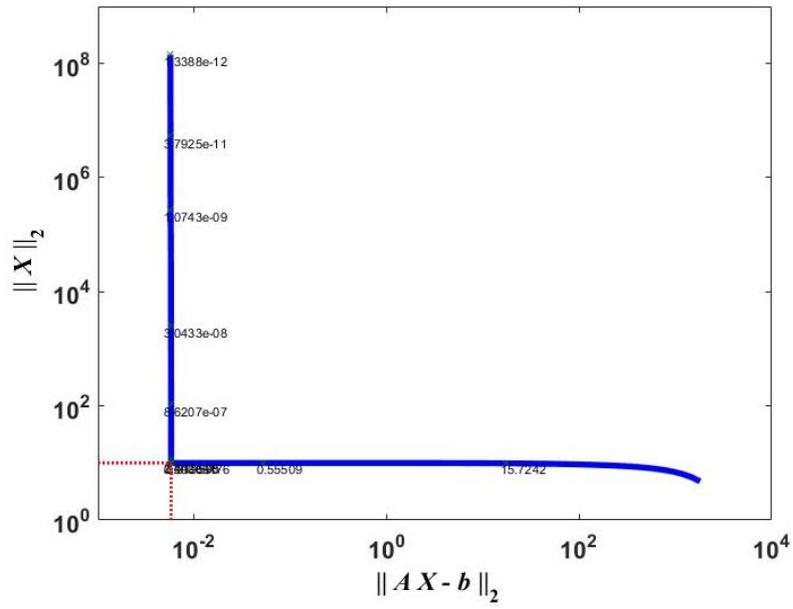

Fig. 2 Regularization Parameter Choice ( L-Curve )

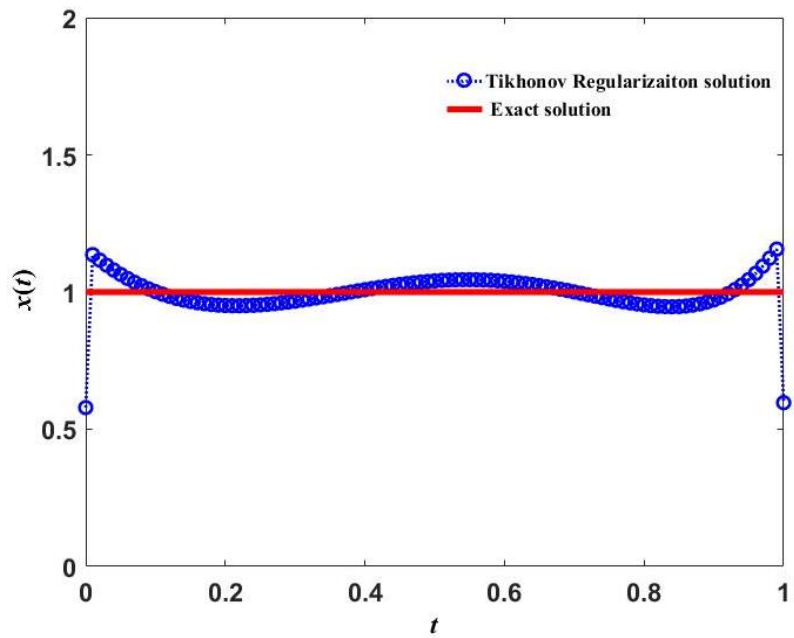

Fig. 3 The results of Tikhonov regularization method

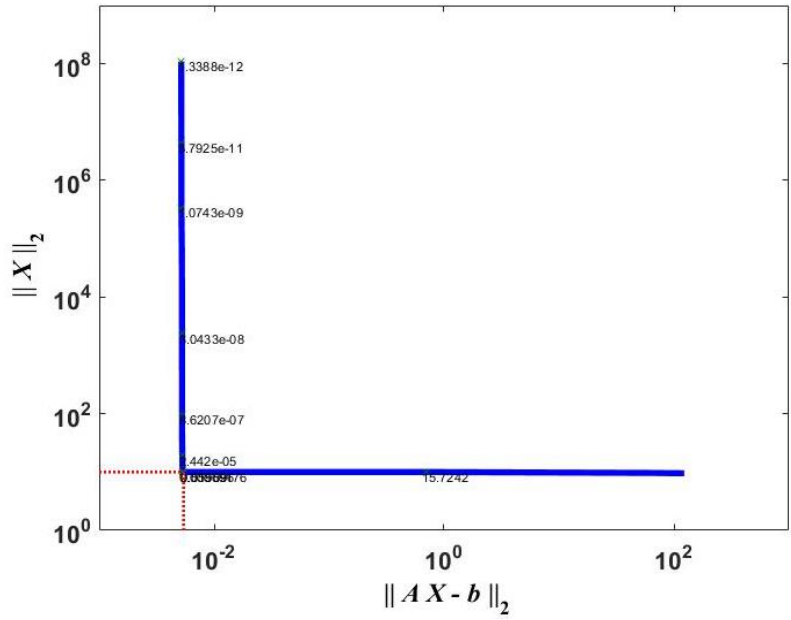

Fig. 4 Regularization Parameter Choice ( L-Curve )

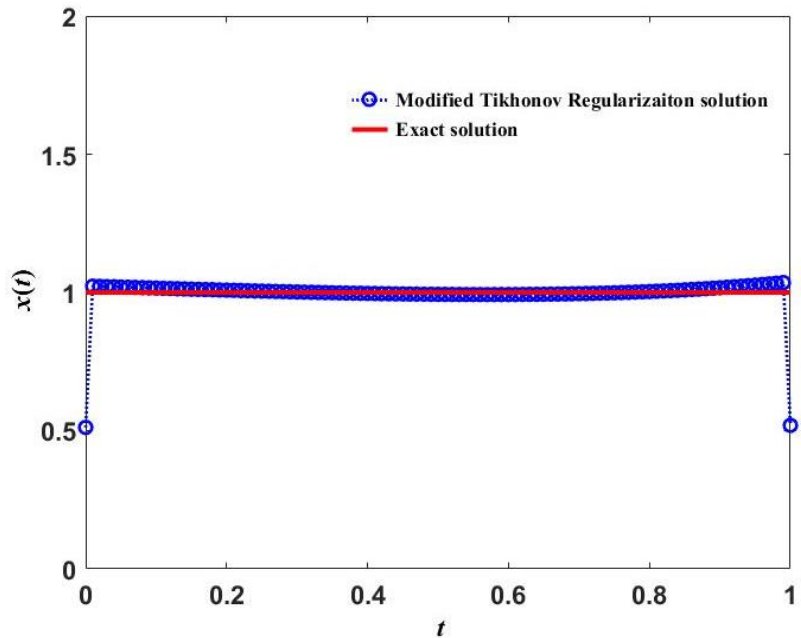

Fig. 5 The results of modified Tikhonov regularization

**Example 4.2** The integral equation of first kind

$$\int_0^1 e^{ts} x(s) \mathrm{d}s = y(t), \quad 0 \leq t \leq 1 \tag{35}$$

with the right-hand side and the solution given by $y(t) = \dfrac{e^{t+2} - 1}{t+2}$ and $x(t) = e^{2s}$.

The parameters involved in the numerical simulation are the same as the previous example.

We also solve this equation with the Tikhonov regularization method and the modified Tikhonov regularization method, respectively. The regularization parameter is chosen by using the L-curve method, and the corresponding curve is shown in Fig.6 and Fig.8. Fig. 7 depicts the comparison between the results obtained by the classical Tikhonov regularization method and exact solution, meanwhile Fig.9 gives the comparison between the results obtained by the modified Tikhonov regularization method and exact solution. From Fig.7 and Fig.9, we can conclude that the modified Tikhonov regularization method work better for this problem.

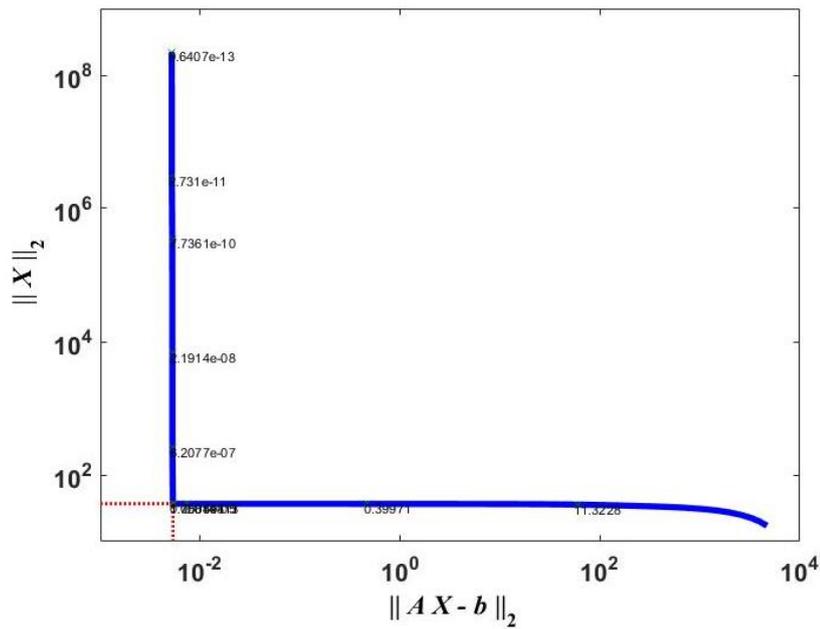

Fig. 6 Regularization Parameter Choice ( L-Curve )

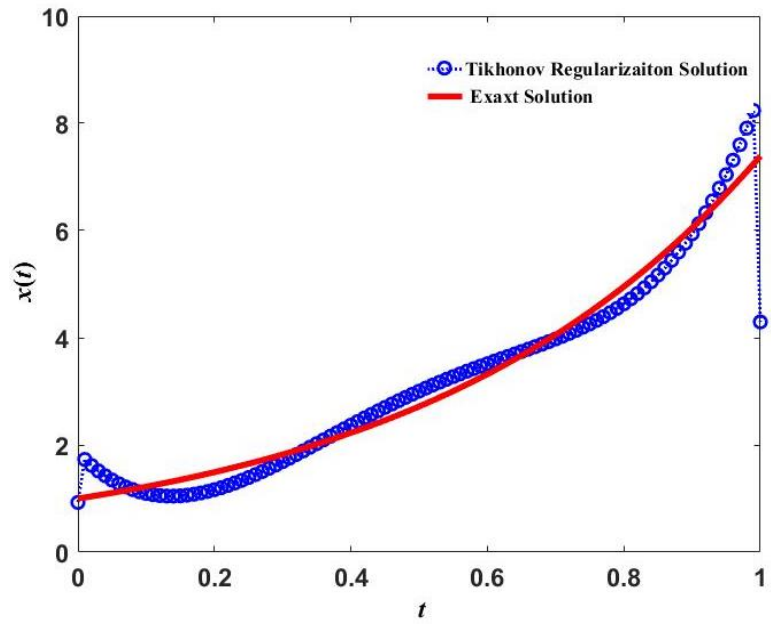

Fig.7 The results of Tikhonov regularization method

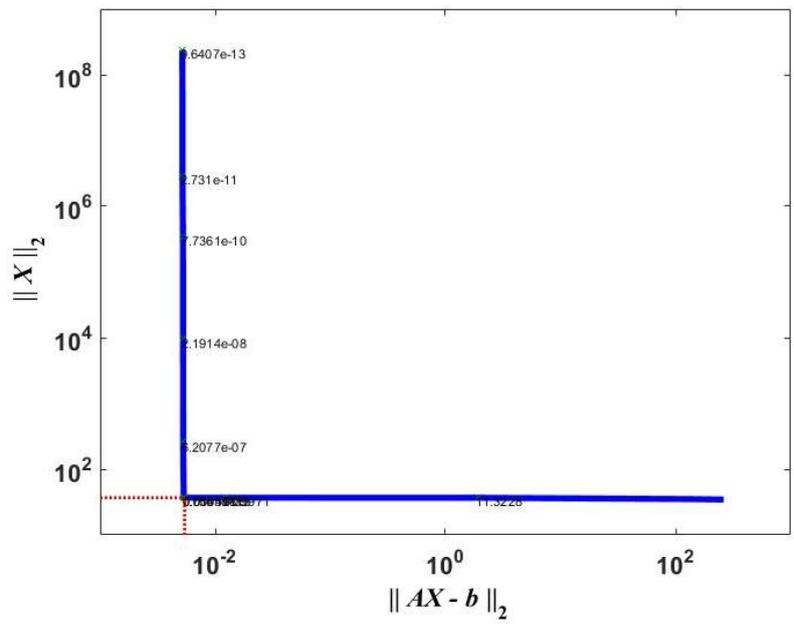

Fig. 8 Regularization Parameter Choice ( L-Curve )

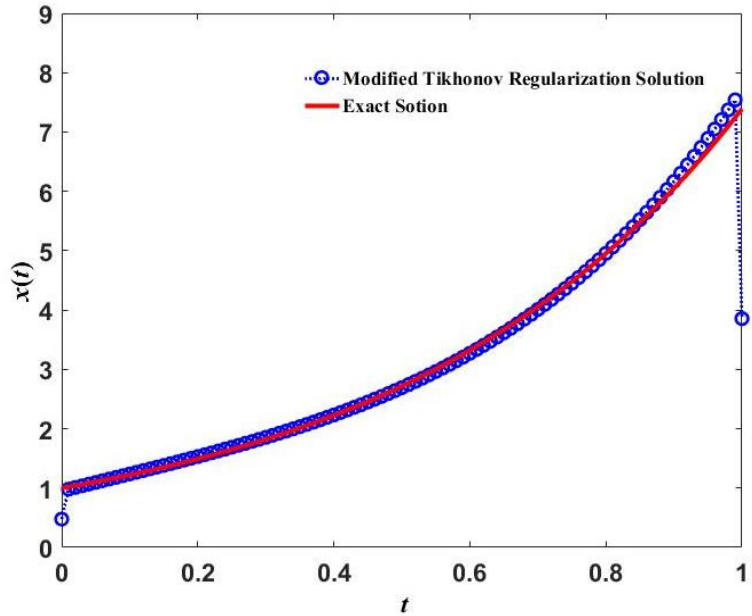

Fig.9 The results of Tikhonov regularization method

## 5. Optimization of drug controlled release from multi-laminated devices

**5.1 The optimization of initial drug concentration**

The initial drug concentration is an essential parameter in the multi-laminated controlled release system, which can affect the drug release greatly. A reasonable initial drug concentration can lead to the drug release with the desired flux. Based on the inverse problem frame (Eq. (6)-Eq.(10) ), the initial drug concentration can be inverted from the known drug release flux. In the following, we will determine the initial drug concentration for three different cases with the Tikhonov regularization method (TRM) and the modified Tikhonov regularization method (MTRM), respectively. The three different desired release profiles are shown in Fig.10.

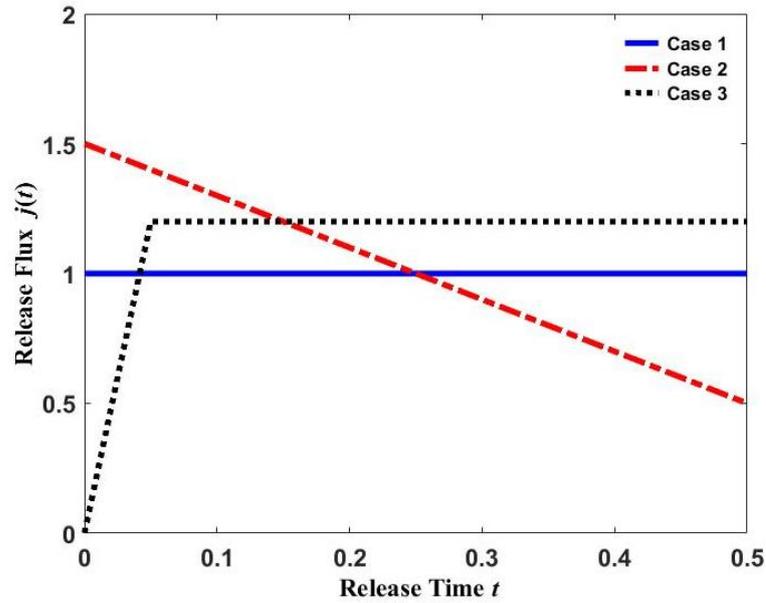

Fig. 10 Three different desired flux. Case 1: Constant release rate. Case 2: Linearly decreasing release rate. Case 3: Nonlinear release rate

**Case 1**

A typical case where the desired flux is assumed to be constant ( $j(t)=1, 0 \leq t \leq 0.5$ ). For the Fredholm integral equation of the first kind (13), set the right-hand side equal to 1, that is, $j(t)=1$. The Tikhonov regularization method and modified Tikhonov regularization method are applied to solve this ill-posed problem. The inverse results are shown in Fig.11, and the computational drug release flux based on the inversed initial drug concentration are depicted in Fig. 12. It is seen that from Fig.12, the optimal release profile obtained using TRM remains flatter pattern at initial stage and has a smaller deviation from ideal case. However, bigger error will appear as the time increase. And the mean square deviation for TRM reaches 0.2174. We can also observe that, for MTRM, although the optimized release fluctuates at the initial stage, the deviation from the desired release remains smaller over the entire computation time. Compared with TRM, the mean square deviation for MTRM is also relatively less, which is 0.1692.

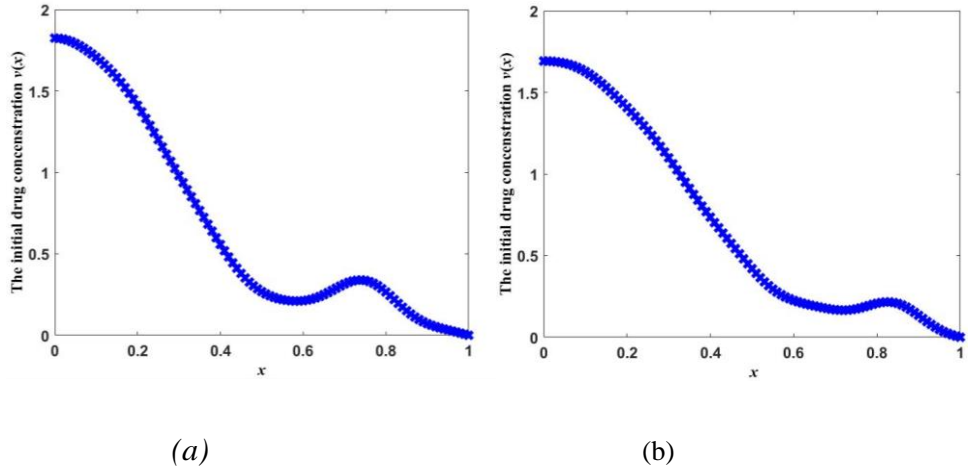

        *(a)*                  (b)

Fig.11 The inverse results ( *a*: Results obtained by TRM, b: Results obtained by MTRM)

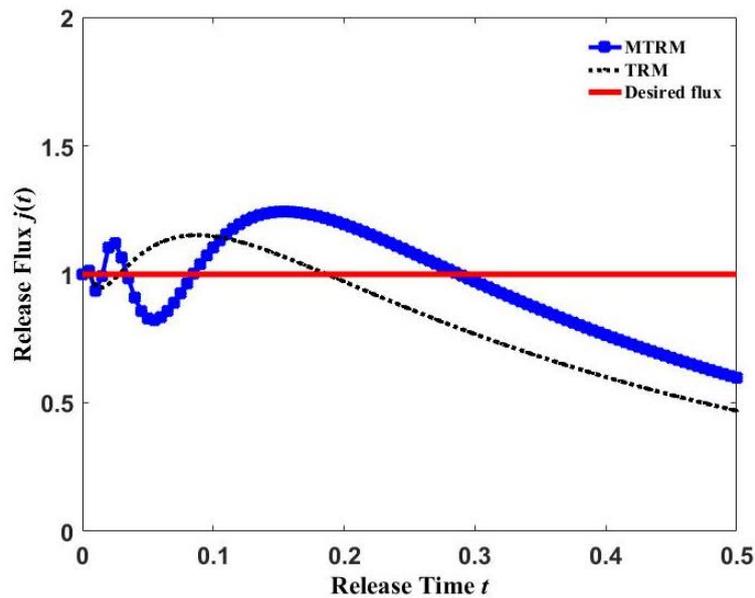

Fig.12 The computational drug release flux based on the inversed initial drug concentration

**Case 2**

    An approximately linearly increasing profile may be desired in some cases, e.g. to build up a tolerance for the chemical being delivered. In the following, the ideal release profile is given with the function $j(t) = 1.5 - 2t$, $0 \leq t \leq 0.5$. The inversed results obtained by using these two different regularization methods are shown in Fig.13. Meanwhile, the optimized release profiles are depicted in Fig. 14. We can conclude that MTRM has obvious advantages over the TRM. In Fig.14, the optimized release profile almost coincides with the desired release flux for MTRM, whereas the performance of TRM remains very poor. The mean square deviations of two methods

are 0.1285 and 0.0164, respectively.

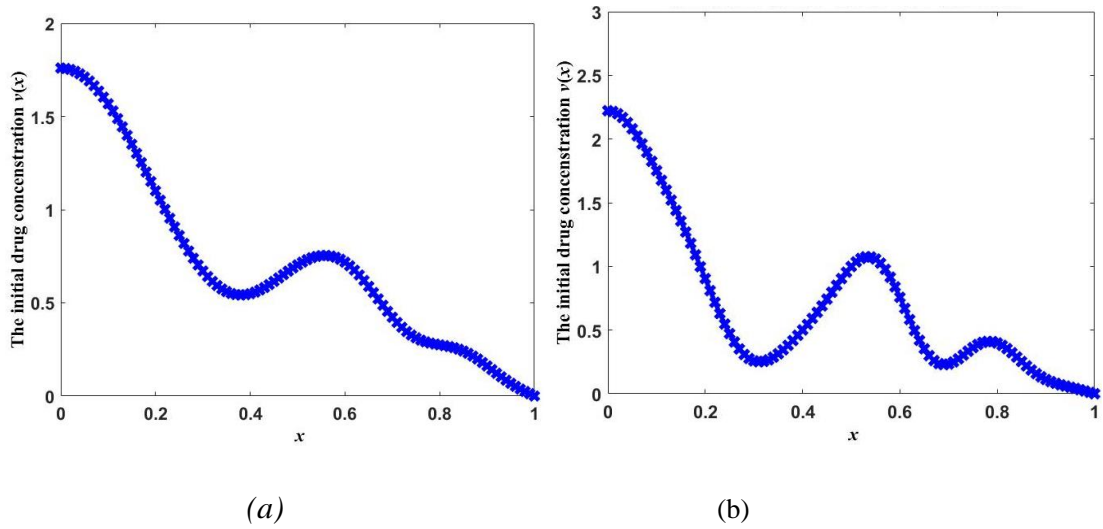

*(a)*                          (b)

Fig.13 The inverse results ( *a*: Results obtained by TRM, b: Results obtained by MTRM)

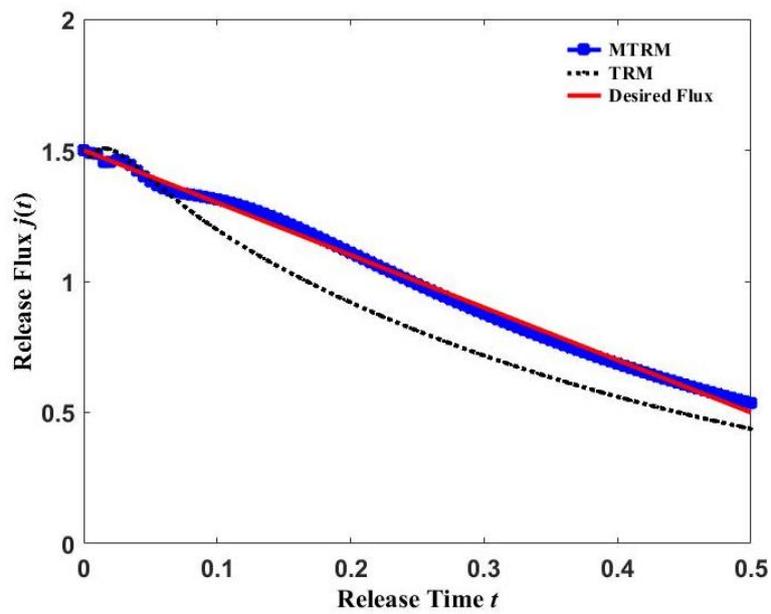

Fig.14 The optimized drug release flux based on the inversed initial drug concentration

**Case 3**

Some situations demand a nonlinearly release rate, e.g., linearly increasing follow by a constant release, and without burst. A typical example is the delivery of the some anticancer drug.

For this case, the release rate function is $j(t) = \begin{cases} 24t & 0 \le t \le 0.05 \\ 1.2 & 0.05 < t \le 0.5 \end{cases}$. The inversed results obtained by using these two different regularization methods are shown in Fig.15. Meanwhile, the optimized release profiles are depicted in Fig. 16. From these two figures, we can see that the performance of MTRM is better than that of TRM. The mean square deviations of two methods are 0.3019 and 0.2164, respectively.

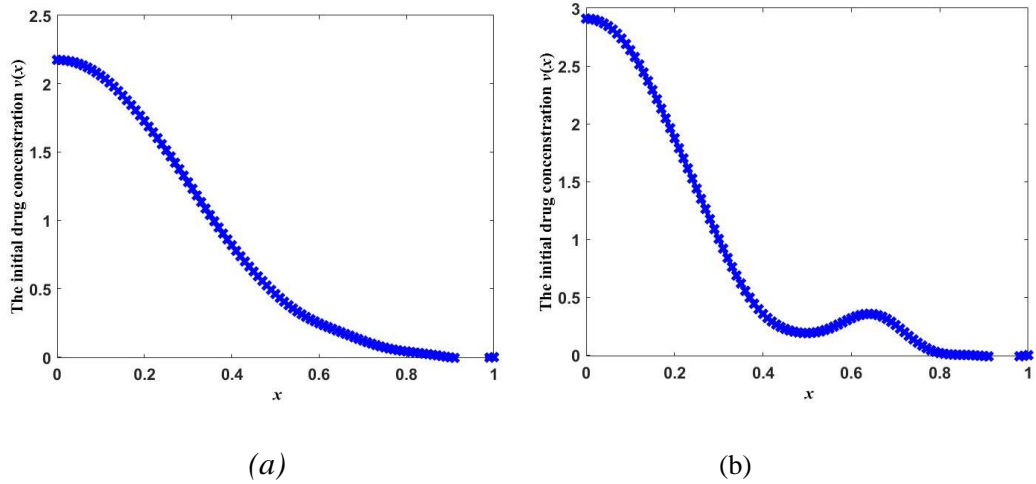

(a)              (b)

Fig.15 The inverse results ( *a*: Results obtained by TRM, b: Results obtained by MTRM)

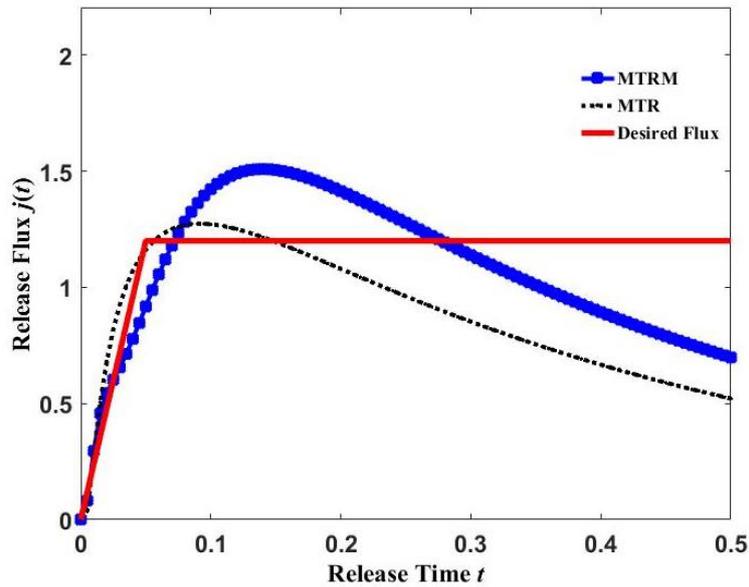

Fig.16 The optimized drug release flux based on the inversed initial drug concentration

## 5.2 Anti-noise property analysis

In practice, the right-hand side $j(t)$ of the Fredholm integral equation of the first kind (13) is never known exactly but only up to an error, say, $\delta > 0$. Therefore, we will consider the inverse problem with a perturbed data.

Assume that we know $\delta > 0$ and with $j^\delta(t)$ satisfy $\left|j^\delta(t) - j(t)\right| \leq \delta$. It is our aim to solve the perturbed equation

$$2\int_0^1 \left[\sum_{k=0}^\infty (-1)^k \left(k + \frac{1}{2}\right)\pi e^{-\left[\left(k+\frac{1}{2}\right)\pi\right]^2 t} \cos\left(k + \frac{1}{2}\right)\pi x\right] v^\delta(x) \mathrm{d}x = j^\delta(t), \qquad (36)$$

We choose $\delta = 0.1$. The inversed results based on the equation (19) by using MTRM for three different cases are shown in Fig.17-Fig.19. These three pictures demonstrate that the modified Tikhonov regularization method we proposed can still work well. It means that the MTRM has a good anti-noise property.

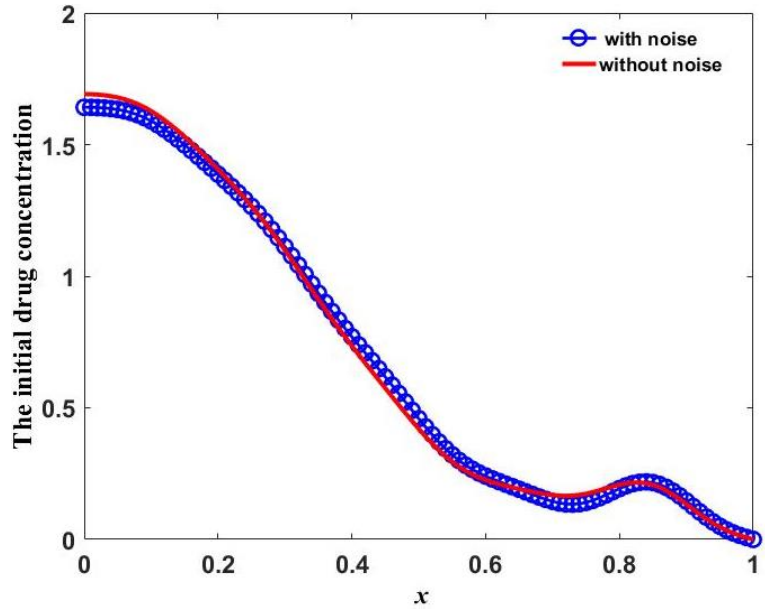

Fig. 17 The inverse initial drug concentration for Case 1 with $\delta = 0.1$

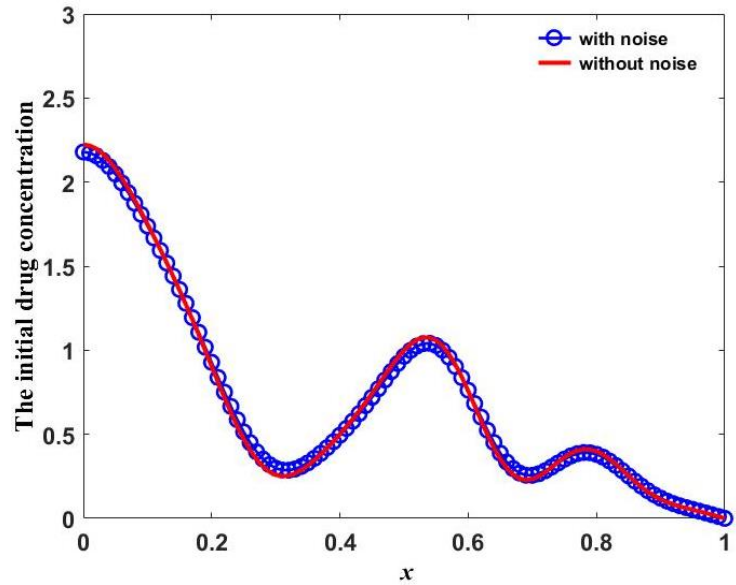

Fig. 18 The inverse initial drug concentration for Case 2 with $\delta = 0.1$

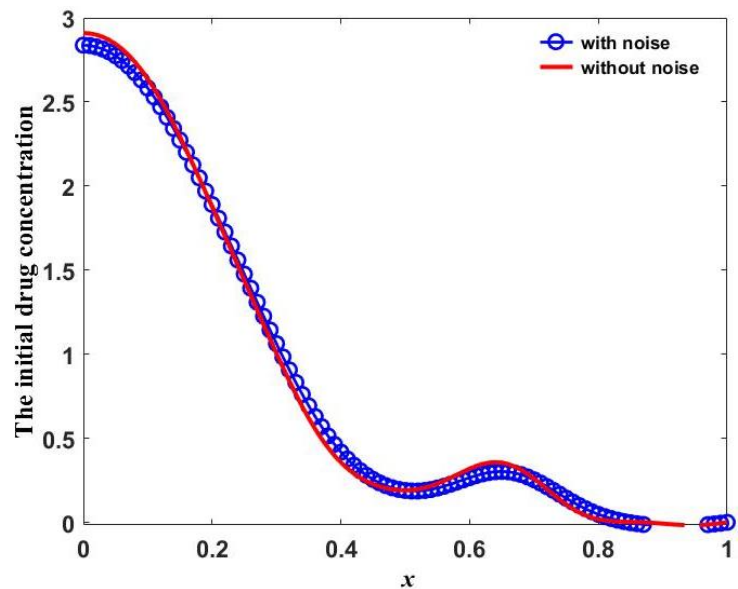

Fig. 18 The inverse initial drug concentration for Case 3 with $\delta = 0.1$

## 6. Conclusion

We have proposed a new viewpoint to solve the optimization problem of drug release based on the multi-laminated drug controlled release device, that is, inverse problem frame. The objective of this paper is to show that the inverse problem frame is effective for the optimization problem of drug release. Based on the inverse problem frame, the optimization problem of drug release can be transformed into the solution of diffusion equation initial value inverse problem, and then be converted to the solution of the Fredholm integral equation of first kind. This is an ill-posed problem, which need the special regularization method to solve.

To solve this ill-posed problem, we presented a modified Tikhonov regularization method is proposed by constructing a new regularizing filter based on the singular value theory of compact operator and proved the convergence and the optimal asymptotic order of the regularized solution. Furthermore, for three various desired release flux, the modified Tikhonov regularization method is applied to inverse the initial drug concentration. As seen in the examples, the method proposed in this paper has been successful at inverting the initial drug concentration. This demonstrates that the modified Tikhonov regularization approach is well suited to solving this ill-posed nonlinear problem.

Also shown in this paper is the result that the modified Tikhonov regularization method has a better anti-noise property for the initial drug concentration estimation. With10% noise, the results obtained with the MTRM are satisfactory.

In all, the better results obtained mean that the solution frame based on the inverse problem for the optimization problem of drug release based on the multi-laminated drug controlled release device exhibits its effectiveness and superiority to some extent in both theoretical research and numerical simulation. There is good potential that the method can be employed to solve more complicated cases, such as multi-parameter identification and high dimensional cases. And this is an important direction for us to face in future.

## Acknowledgments

The work was supported by the National Science Foundation of China, under Grant No. 41004052 and supported by the Natural Science Foundation of Guangdong Province, under Grant No. 2017A030313280.